\documentclass{amsart}

\theoremstyle{definition}

\theoremstyle{remark}

\numberwithin{equation}{section}

\reversemarginpar

\newcommand{\nn}{\nonumber}
\newcommand{\no}{\noindent}

\newcommand{\realpart}{\mathop{\rm Re}\nolimits}

\newcommand{\lp}{\ln \sqrt{2 \pi}}
\newcommand{\ba}{\begin{eqnarray}}
\newcommand{\ea}{\end{eqnarray}}
\newcommand{\bml}{\begin{multline}}
\newcommand{\eml}{\end{multline}}

\newcommand{\ift}{\int_{0}^{\infty}}
\newcommand{\ione}{\int_{0}^{1}}
\newcommand{\ionehalf}{\int_{0}^{1/2}}

\newcommand{\allR}{\mathbb{R}}
\newcommand{\allN}{\mathbb{N}}
\newcommand{\nnN}{\mathbb{N}_{0}}

\newcommand{\stnk}[2]{\left\{\begin{matrix} #1\\#2\end{matrix}\right\}}
\newcommand{\st}{{ \quad \quad \quad  }}
\newcommand{\stoli}{}
\newcommand{\stboth}{}

\begin{document}

\begin{center}
Submitted to {\bf The Ramanujan Journal}. July 2001.
\end{center}

\title[Hurwitz zeta function] {On some integrals involving the Hurwitz
zeta function: part 2}

\author{Olivier Espinosa}
\address{Departamento de F\'{\i}sica,
Universidad T\'{e}cnica Federico Santa Mar\'{\i}a, Valpara\'{\i}so, Chile}
\email{espinosa@fis.utfsm.cl}

\author{Victor H. Moll}
\address{Department of Mathematics,
Tulane University, New Orleans, LA 70118}
\email{vhm@math.tulane.edu}

\subjclass{Primary 33}

\date{\today}


\begin{abstract}
We establish a series of indefinite integral formulae involving
the Hurwitz zeta function and other elementary and special
functions related to it, such as the Bernoulli polynomials, $\ln
\sin (\pi q)$, $\ln \Gamma(q)$ and the polygamma functions. Many of
the results are most conveniently formulated in terms of a family
of functions $A_k(q):=k\zeta'(1-k,q)$, $k\in\allN$, and a family of
polygamma functions of negative order, whose properties we study
in some detail.
\end{abstract}

\maketitle


\newtheorem{Definition}{\bf Definition}[section]
\newtheorem{Thm}[Definition]{\bf Theorem}
\newtheorem{Lem}[Definition]{\bf Lemma}
\newtheorem{Cor}[Definition]{\bf Corollary}
\newtheorem{Prop}[Definition]{\bf Proposition}
\newtheorem{Example}[Definition]{\bf Example}

\section{Introduction} \label{S:intro}

The Hurwitz zeta function, defined by
\ba
\zeta(z,q) & = & \sum_{n=0}^{\infty} \frac{1}{(n+q)^{z}}
\ea
\no
for $z \in \mathbb{C}, \, \realpart{z} > 1$ and $q \neq 0, -1, -2, \cdots,$
admits the integral representation
\ba
\zeta(z,q) & = & \frac{1}{\Gamma(z)} \int_{0}^{\infty}
\frac{e^{-qt}}{1-e^{-t}} t^{z-1} dt,
\ea
\no
where $\Gamma(z)$ is Euler's gamma function, which is valid for
$\realpart z > 1$ and $\realpart q >0$, and can be used to
prove that $\zeta(z,q)$
has an analytic extension to the whole complex plane except for a
simple pole at $z=1$.

For $\realpart z < 0$, $\zeta(z,q)$ admits the following Fourier
representation, originally derived by
Hurwitz, in the range $0 \le q \le 1$:
\begin{multline}
\label{fourier}
\zeta(z,q) = \frac{2 \Gamma(1-z)}{(2 \pi)^{1-z}}\;
\left[ \sin \left( \frac{\pi z}{2} \right) \sum_{n=1}^{\infty}
\frac{\cos( 2 \pi q n)}{n^{1-z}} +
\cos \left( \frac{\pi z}{2} \right)
\sum_{n=1}^{\infty} \frac{ \sin( 2 \pi q n) }{ n^{1-z} }
\right].
\end{multline}
\no
This representation was used in \cite{esmo1} to obtain
several definite integral formulae involving $\zeta(z,q)$. A
derivation of (\ref{fourier}) can be found in \cite{ww}, page 268.
An alternative proof, based upon the representation
\ba
\zeta(z,q) & = & \frac{q^{1-z}}{z-1} + \frac{q^{-z}}{2} - z \ift
\frac{\{ t \} - \tfrac{1}{2} }{(t+q)^{z+1}} \, dt,
\label{repre1}
\ea
\no
where $\{ t \}$ is the fractional part of $t$, has been given by
Berndt \cite{hurber}.
The expression (\ref{repre1}) is
employed in \cite{hurber} to give short proofs of several classical
formulae, including Lerch's beautiful expression
\ba
\ln \Gamma(q) & = & \zeta'(0,q) - \zeta'(0). \label{lerch}
\ea

In this paper we continue the work, initiated in \cite{esmo1}, on the explicit
evaluation of integrals involving $\zeta(z,q)$.
Special cases of $\zeta(z,q)$ include the  Bernoulli polynomials,
\ba
\label{zetaber}
B_{m}(q) & = & -m \, \zeta(1-m,q),  \quad m \in \mathbb{N},
\ea
\no
defined by their generating function
\ba
\frac{te^{qt}}{e^{t}-1} & = & \sum_{m=0}^{\infty} B_{m}(q)
\frac{t^{m}}{m!},
\label{generating1}
\ea
\no
and given explicitly in terms of the Bernoulli numbers $B_{k}$ by
\ba \label{berpoly}
B_{m}(q) & = & \sum_{k=0}^{m} \binom{m}{k} B_{k}q^{m-k};
\ea
the digamma function,
\ba\label{digamma-hurwitz}
\psi(q):=\frac{d}{dq} \ln\Gamma(q)=\lim_{z\to 1}
\left[\frac{1}{z-1}-\zeta(z,q)\right];
\ea
and the polygamma functions,
\ba
\psi^{(m)}(q)= (-1)^{m+1}\,{m!}\,\zeta(m+1,q),\quad
m\in\allN,
\label{polygamma-def}
\ea
defined by
\ba \psi^{(m)}(q) & := & \frac{d^{\,m}}{dq^{\,m}} \psi(q),\quad
m\in\allN.
\label{polyga}
\ea

\medskip

An important property of the Hurwitz zeta function, which will be
essential for the indefinite integral evaluations presented in
Section \ref{indefinite}, is the following:
\ba\label{derivative}
\frac{\partial}{\partial q}\, \zeta(z,q) =  -z \,\zeta(z+1,q).
\ea

\medskip

The rest of this paper is organized as follows. In Section
\ref{indefinite} we consider the evaluation of indefinite
integrals of functions of the form $f(q) \zeta(z,a+bq)$, using a
simple integration by parts approach. In Section
\ref{speciala} we introduce and study some of the properties of
two families of functions related to the first derivative with
respect to the argument $z$ of the Hurwitz zeta function
$\zeta(z,q)$, evaluated at $z$ equal to nonpositive integers.
These functions appear in connection to the indefinite
integrals involving polygamma and negapolygamma functions, as well
as $\ln\Gamma(q)$ and $\ln\sin(\pi q)$, considered in
Section \ref{polygamma}. These, in turn, are derived from the
formulae obtained in Section \ref{indefinite} either by direct
differentiation or by taking the appropriate limits. Finally, in
Section \ref{definite} we use some of the indefinite integral
formulae to rederive some of the definite integral evaluations
obtained in reference \cite{esmo1} and to present some new
analogous formulae.
\medskip

\section{The evaluation of indefinite integrals}
\label{indefinite}

In this section we discuss a  method to evaluate primitives of
functions of the form $f(q) \zeta(z,a+bq)$.  This is illustrated
in the cases where $f$ is a polynomial and an exponential
function. The resulting evaluations can be taken as a starting
point to derive similar formulae involving other special functions
in place of $\zeta(z,a+bq)$. For instance,
differentiation with respect to the parameter $z$ leads,
in view of Lerch's result (\ref{lerch}), to the evaluation of
primitives involving the weight $\ln \Gamma(a+bq)$, and thus also
$\ln \sin\pi q$, by virtue of the reflection formula for the gamma
function. Also, the limit $z\to m\in\allN$ leads to the
evaluation of primitives involving the polygamma function
$\psi^{(m-1)}(q)$. We shall present these results, from a slightly
more general point of view, in Section \ref{polygamma}.
\medskip

\begin{Thm} \label{master0}
Let $r \in \allN, \; f$ be $r$-times differentiable and
$a, \, b \in \mathbb{R}$. Then
\ba
\int f(q) \zeta(z,a+bq) \, dq & = & \sum_{k=1}^{r} (-1)^{k+1}
\frac{f^{(k-1)}(q) \, \zeta(z-k,a+bq)}{b^{k} \, (1-z)_{k}}   \label{master}
\\
 & + & \frac{(-1)^{r} }{b^{r} (1-z)_{r}} \int f^{(r)}(q) \,
\zeta(z-r,a+bq) \, dq.  \nn
\ea
\end{Thm}
\begin{proof}
Observe that
\ba
\frac{\partial}{\partial q} \zeta(z-1,a+bq) & = & b(1-z) \zeta(z,a+bq),
\ea
\no
so that integration by  parts yields
\ba
\int f(q) \zeta(z,a+bq) \, dq & = & \frac{f(q) \zeta(z-1,a+bq)}{b(1-z)}
\nn \\
& - & \frac{1}{b(1-z)} \int f'(q) \, \zeta(z-1,a+bq) \, dq. \nn
\ea
\no
The expression (\ref{master}) follows by repeating this procedure.
\end{proof}

\medskip

We now produce the evaluation of certain indefinite integrals by
choosing appropriate
functions $f$ in Theorem \ref{master0}.  \\

\begin{Example}
\label{example22}
Let $n \in \nnN$ and $a, \, b \in \allR$. Then
the moments of $\zeta(z,q)$ are given by
\begin{multline}
\label{mom1}
\int q^{n} \zeta(z,a+bq) \, dq =
n! \sum_{j=0}^{n} \frac{(-1)^j q^{n-j} }{b^{j+1}(1-z)_{j+1} (n-j)!}
\zeta(z-j-1,a+bq).
\end{multline}
\end{Example}
\begin{proof}
The case $n=0$ is simply the known result
\ba\label{Hurwitz-primitive}
\int \zeta(z,a+bq) \, dq & = & \frac{\zeta(z-1,a+bq)}{b(1-z)}.
\ea
For $n\ge 1$, the function $f(q) = q^{n}$ satisfies $f^{(k-1)}(q) =
n!q^{n-k+1}/(n-k+1)!$, for $k\le n$.
Then (\ref{master}), with $r=n$,  yields
\ba
\int q^{n} \zeta(z,a+bq) \, dq & = & \sum_{k=1}^{n}
\frac{(-1)^{k+1} n! q^{n-k+1}}{(n-k+1)! b^{k} (1-z)_{k} } \zeta(z-k,a+bq)
\st \nn \\
 & + & \frac{(-1)^{n} \, n!}{b^{n} (1-z)_{n}} \int \zeta(z-n, a+bq) \, dq,
\st \nn
\ea
\no
so that (\ref{mom1}) follows from (\ref{Hurwitz-primitive}).
\end{proof}
\no
In a similar fashion we obtain:
\begin{Example}
Let $m \in \nnN$ and $a, \, b, \, c, \, d\in \allR$. Then
\begin{equation}
\label{momber1}
\int B_m(c+dq) \zeta(z,a+bq) \, dq =
m! \sum_{j=0}^{m} \frac{(-1)^j d^{\,j} B_{m-j}(c+dq) }{b^{j+1}(1-z)_{j+1} (m-j)!}
\zeta(z-j-1,a+bq).
\end{equation}
\end{Example}
\begin{proof}
Same as the proof for Example $3.2$ with
\[
\frac{{d^{\,k - 1} }}{{dq^{\,k - 1} }}B_m (c+dq) =
\frac{{m! \, d^{k-1} }}{{(m - k + 1)!}}B_{m - k + 1} (c+dq).
\]
\end{proof}
\medskip

\no
{\bf Definition}.  The family of functions
$\mathfrak{F} := \{ f_{j}(q): \, j \in \mathbb{N}
\}$ is said to be {\em closed under primitives} if for each $j$
 the primitive of
$f_{j}(q)$ can be written as a {\em finite} linear combination of the
elements of $\mathfrak{F}$.  Naturally, the family $\mathfrak{F}$ is allowed
to depend on a finite number of parameters, as in the
next example.  \\

\no
\begin{Example}
Example \ref{example22} shows that
\begin{multline}
{\mathfrak{F}}_{a,b}  :=  \{ P_{j}(q) \zeta(z-m,a+bq): j, \, m \in
\mathbb{N} \text{ and } P_{j} \text{ is  a
 polynomial in } q \\
\text{ of degree } j\, \} \nn
\end{multline}
\no
is closed under primitives.  This follows from
\begin{multline}
\label{mom2}
\int q^{n} \zeta(z-m,a+bq) \, dq =
n! \sum_{j=0}^{n} \frac{(-1)^j q^{n-j}
\zeta(z-m-1-j,a+bq)}{b^{j+1}(m+1-z)_{j+1} (n-j)!},
\end{multline}
which is a variation of (\ref{mom1}).
\end{Example}

\medskip

\begin{Example}
The moments of the Bernoulli polynomials are given by
\begin{multline}
\label{mom3}
\int q^{n} B_{m}(a+bq) \, dq =
\frac{n!m!}{(n+m+1)!} \sum_{j=0}^{n}
\frac{(-1)^{j} q^{n-j}}{b^{j+1}} \binom{m+n+1}{n-j}
B_{m+j+1}(a+bq).
\end{multline}
\end{Example}
\begin{proof}
Use the identity (\ref{zetaber}) in (\ref{mom1}).
\end{proof}

\medskip

\begin{Example}
\label{intodd}
Let $n \in \mathbb{N} $ be odd. Then
\begin{multline}
\int \zeta(z-n,q) \zeta(z,q) dq =
\frac{1}{2}
\sum_{k=1}^{n} \frac{(z-n)_{k-1}}{(1-z)_{k}} \zeta(z-k,q)
\zeta(z-n+k-1,q). \st
\end{multline}
\end{Example}

\begin{proof}
The Hurwitz zeta function satisfies
\ba
\label{hur100}
\frac{\partial^{k-1}}{\partial q^{k-1}} \zeta(z-n,q)=
(-1)^{k-1} (z-n)_{k-1} \zeta(z-n+k-1,q),
\ea
\no
so the result follows from (\ref{master}) with $r=n$, since in
that case the integral on the right-hand side equals the one on
the left-hand side, except for the prefactor
$(z-n)_n/(1-z)_n=(-1)^n$. \\
\end{proof}
\no
{\bf Note.} In view of the identity
\[
\frac{{(z - n)_{k - 1} }}{{(1 - z)_k }} = ( - 1)^{n + 1} \frac{{(z
- n)_{n - k} }}{{(1 - z)_{n - k + 1} }},
\]
it is easily seen that the terms in the sum on the right-hand side
of (\ref{hur100}) are equal in pairs, except for the central term
$k=r$, where $r\in\allN$ is defined by $n=2r-1$.
Therefore we have the alternative formula
\begin{multline}
\int \zeta(z-2r+1,q) \zeta(z,q) dq =
\frac{(z-2r+1)_{r-1}}{2\,(1-z)_{r}} \zeta^2(z-r,q)
\\
+\sum_{k=1}^{r-1} \frac{(z-2r+1)_{k-1}}{(1-z)_{k}} \zeta(z-k,q)
\zeta(z-(2r-k),q).
\end{multline}

\medskip

%

\no
{\bf Note}. We have been unable to evaluate the integral in Example
\ref{intodd} for
the case $n$ even. Thus the question of whether the family
\ba
\mathfrak{F}_{z} & := & \{ \zeta(z-n,q) \zeta(z-m,q): \, n,m \in
\mathbb{N} \}
\ea
\no
is closed under primitives remains to be decided.  \\

\medskip

\begin{Example}
Let $a, \, b \in \mathbb{R}$. Then
\ba
 \int e^{q} \zeta(z,a+bq) \, dq & = & e^{q}
\sum_{j=0}^{\infty} \frac{(-1)^{j}}{b^{j+1} (1-z)_{j+1} } \;
\zeta(z-1-j,a+bq). \label{exp1}
\ea
\end{Example}
\begin{proof}
Divide (\ref{mom1}) by $n!$ and then sum over $n$ to produce
\ba
\int e^{q} \zeta(z,a+bq) \, dq & = &
\sum_{n=0}^{\infty} \sum_{j=0}^{n} \frac{(-1)^{j} q^{n-j} }{b^{j+1} (1-z)_{j+1}
(n-j)!} \, \zeta(z-1-j,a+bq). \nn
\ea
\no
The result follows by interchanging the order of summation.
\end{proof}

\noindent
{\bf Note}. We have been unable to produce a {\em finite expression} for
the integral in (\ref{exp1}).
\medskip

\begin{Example}
Let $m \in \mathbb{N}$. Then
\ba
\int e^{q} B_{m}(a+bq) \, dq & = & m! e^{q} (-1)^m
\sum_{j=0}^{m} \frac{(-1)^{j}}{j!} b^{m-j} B_{j}(a+bq).
\ea
\end{Example}
\begin{proof}
Use the identity (\ref{zetaber})  and $(m)_{j+1} = (m+j)!/(m-1)!$
in (\ref{exp1}) to produce
\ba
\int e^{q} B_{m}(a+bq) \, dq  & = & m!e^{q}(-1)^{m+1} \sum_{j=m+1}^{\infty}
\frac{(-1)^{j}}{j!} b^{m-j} B_{j}(a+bq).
\ea
\no
The generating function (\ref{generating1}) is now employed to see that the
sum from $j=0$ to infinity is independent of $q$, so it is absorbed into the
implicit constant of integration.
\end{proof}

\medskip

\section{The function $A_{k}(q)$ and negapolygammas}
\label{speciala}

In this section we consider the function
\ba
\label{defAk}
A_{k}(q) := k\,\frac{\partial}{\partial z}\, \zeta(z,q)\Big{|}_{z=1-k}
\ea
\no
for $k \in \allN$. The function $A_1(q)$ has a simple explicit
form,
\\
\ba\label{def-A1}
A_1(q)=\zeta'(0,q)=\ln\Gamma(q)+\zeta'(0),
\ea
\\
in view of Lerch's result (\ref{lerch}).
These functions appear in all of the formulae
for the indefinite integrals involving the loggamma
and the logsine functions studied in Section \ref{polygamma}.
The derivative of the Hurwitz zeta function has appeared before in
connection to integrals of $\ln\Gamma(q)$ \cite{gosper}, and in a
number of related contexts, such as the studies of polygamma
functions of negative order \cite{adam1}, the Barnes function
\cite{adam2} and the multiple gamma function \cite{vardi},
and other unrelated ones such as the evaluation of sums of
the type $\sum_{m\ge 2} z^m\zeta(m,\alpha)$ \cite{kanemitsu}.
Here, we will rediscover, from a more general point of view, the
intimate relationship existing between the functions $A_k(q)$ and
polygamma functions of negative order, first pointed out in
reference \cite{gosper}.
\medskip

The nonelementary behavior of $A_k(q)$ can be considered to be
contained in the range $0\le q\le 1$, since
for $q>1$ the value of $A_k(q)$ can be obtained by repeated use
of the following result:
\begin{Lem}
The function $A_{k}(q)$ satisfies
\ba
A_{k}(q+1) & = & A_{k}(q)+ k\, q^{k-1} \ln q. \label{recu2}
\ea
\end{Lem}
\begin{proof}
Differentiate both sides of the identity
\ba
\label{zeta-shift}
\zeta(z,q)=\frac{1}{q^z}+\zeta(z,q+1)
\ea
with respect to $z$ at $z=1-k$.
\end{proof}
As an immediate consequence of  property (\ref{recu2}) and
definition (\ref{defAk}) we have
\begin{Lem}
For $k\ge 2$,
\ba
A_{k}(0) = A_{k}(1) = k\,\zeta'(1-k).\label{Ak-at-zero}
\ea
\end{Lem}
\begin{proof}
Take the limit $q\to 0$ in (\ref{recu2}) and use the fact that
$\zeta(z,1)=\zeta(z)$.
\end{proof}
\begin{Lem}
For $k\in\allN$,
\ba
\ione A_k(q)\,dq=0.\label{Ak-int-0-1}
\ea
\end{Lem}
\begin{proof}
This is a direct consequence of
\ba
\ione \zeta(z,q)\,dq=0,
\ea
valid for $\realpart{z}<1$.
\end{proof}

\begin{Lem}\label{der-Ak}
For $k\in\allN$, the functions $A_{k}(q)$ satisfy
\ba
A_{k+1}'(q) & = & (k+1)\left[A_{k}(q)+\frac{1}{k}B_{k}(q)\right]. \label{recu1}
\ea
\end{Lem}
\begin{proof}
We have
\ba
A_{k+1}'(q) & = & (k+1)\frac{\partial}{\partial q} \frac{\partial }{\partial z}
\zeta(z,q)\Big{|}_{z=-k} =
(k+1)\frac{\partial}{\partial z} \frac{\partial }{\partial q}
\zeta(z,q)\Big{|}_{z=-k}  \nn \\
 & = & (k+1)\frac{\partial}{\partial z} \left[ -z \zeta(z+1,q) \right]
\Big{|}_{z=-k} \nn \\
& = & (k+1)\left[A_{k}(q)+\frac{1}{k}B_{k}(q)\right].  \nn
\ea
\end{proof}
\no
The following results can sometimes be used to simplify a formula:
\ba
\zeta'(-2n) & = & (-1)^{n} \frac{(2n)! \, \zeta(2n+1)}{2 \, (2 \pi)^{2n} },
\quad n \in\allN,
\stboth
\label{berzeta2}  \\
\zeta'(0)  & = & -\lp.
\stboth
\label{berzeta3}
\ea
\begin{Lem}
For $k\in\allN$,
\ba
\label{Ak-at-onehalf}
A_{k}(\tfrac{1}{2}) = ( - 1)^{k-1} B_{k}
2^{1-k}\ln 2 - (1 - 2^{1 - k} )k\,\zeta'(1-k).
\ea
In particular,
\ba
A_1(\tfrac{1}{2})&=&-\frac{1}{2}\ln 2, \\
A_{2n+1}(\tfrac{1}{2}) &=&
( - 1)^{n + 1} \frac{{(1 - 2^{ -2n} )(2n+1)!\zeta (2n + 1)}}
{{2(2\pi )^{2n} }}\quad n\ge 1.
\ea
\end{Lem}
\begin{proof}
Differentiate $\zeta(z,\tfrac{1}{2})=\left(2^z-1\right)\zeta(z)$ at $z=1-k$ and
use the identity
\ba
\zeta(1-k) = \frac{(-1)^{k+1} B_{k}}{k}, \quad k \in\allN.
\stboth
\label{berzeta1}
\ea
\end{proof}
\no

Hurwitz's Fourier representation (\ref{fourier}) for $\zeta(z,q)$
in the range $0\le q\le 1$ and negative $z$ is absolutely
convergent and thus can be used to directly obtain a Fourier
representation for the function $A_{k}(q)$ in the range $0\le q\le
1$.

Define
\ba
C(z,q) = \sum_{n=1}^{\infty}\frac{\cos (2 \pi n q) }{n^{z}} & \text{ and } &
S(z,q) = \sum_{n=1}^{\infty}\frac{\sin (2 \pi n q) }{n^{z}}, \label{CandS}
\ea
\no
so that
\ba
\frac{\partial}{\partial z} C(z,q)
= -\sum_{n=1}^{\infty}\frac{\ln n}{n^{z}} \cos (2 \pi n q) & \text{ and } &
\frac{\partial}{\partial z} S(z,q)
= -\sum_{n=1}^{\infty}\frac{\ln n}{n^{z}} \sin (2 \pi n q). \nn
\ea
\no
Now replace $z$ by $1-z$ in the Fourier representation
(\ref{fourier}) and differentiate to produce
\begin{multline}
\zeta'(1-z,q)  =  \frac{2 \Gamma(z)}{(2 \pi)^{z}} \times\\
\bigg{\{} \left[ \Psi(z) \cos \frac{\pi z}{2} + \frac{\pi}{2}
\sin \frac{\pi z}{2} \right] C(z,q)
 + \left[ \Psi(z) \sin \frac{\pi z}{2} -
\frac{\pi}{2} \cos \frac{\pi z}{2} \right] S(z,q)\\
- \cos \frac{ \pi z}{2} \frac{\partial}{\partial z} C(z,q) -
\sin \frac{ \pi z}{2}  \frac{\partial}{\partial z} S(z,q)
\bigg{\}},
\end{multline}
\no
where $\Psi(z) := \ln 2 \pi - \psi(z)$.

For $z$ a  positive integer, the functions $S(z,q)$ and
$C(z,q)$ are related to the Bernoulli polynomials (in view of
(\ref{fourier}) and (\ref{zetaber})) and the Clausen functions. The
latter are defined by
\ba
Cl_{2n}(x) & = & \sum_{k=1}^{\infty} \frac{\sin kx}{k^{2n}}, \; n
\geq 1
\ea
\no
and
\ba
Cl_{2n+1}(x) & = & \sum_{k=1}^{\infty} \frac{\cos kx}{k^{2n+1}},
\; n \geq 0.
\ea
\no
One has
\ba
S(2m+1,q) &=& \sum_{n=1}^{\infty} \frac{\sin( 2 \pi n q)}{n^{2m+1}} =
\frac{(-1)^{m+1} \, (2 \pi )^{2m+1} }{2(2m+1)!} B_{2m+1}(q),
\label{fourier-berodd}\\
C(2m+1,q) &=& \sum_{n=1}^{\infty} \frac{\cos( 2 \pi n q)}{n^{2m+1}} =
Cl_{2m+1}(2 \pi q),
\label{fourier-clauodd}\\
S(2m+2,q) &=& \sum_{n=1}^{\infty} \frac{\sin( 2 \pi n q)}{n^{2m+2}} =
Cl_{2m+2}(2 \pi q),
\label{fourier-claueven}\\
C(2m+2,q) &=& \sum_{n=1}^{\infty} \frac{\cos( 2 \pi n q)}{n^{2m+2}} =
\frac{(-1)^{m} \, (2 \pi )^{2m+2} }{2(2m+2)!} B_{2m+2}(q).
\label{fourier-bereven}
\ea
\no
The value $z = 2m+1$ yields, upon using
$\psi(k+1) = - \gamma  + H_{k}$, the expression
\begin{multline}
\label{amodd}
A_{2m+1}(q) = ( H_{2m} - \gamma - \ln 2 \pi )\,B_{2m+1}(q)\\
+ (-1)^{m} \frac{2 (2m+1)!}{(2 \pi)^{2m+1}}
\left[ \sum_{n=1}^{\infty} \frac{ \ln n}{n^{2m+1}} \sin( 2 \pi n q) +
\frac{\pi}{2} \sum_{n=1}^{\infty} \frac{\cos( 2 \pi n q)}{n^{2 m+1} }
\right].
\end{multline}
\no
Similarly, for $z = 2m+2$ we find
\begin{multline}
\label{ameven}
A_{2m+2}(q) = ( H_{2m+1} - \gamma - \ln 2 \pi )
\,B_{2m+2}(q)\\
+ (-1)^{m+1} \frac{2(2m+2)!}{(2 \pi)^{2m+2}}
\left[ \sum_{n=1}^{\infty} \frac{ \ln n}{n^{2m+2}} \cos( 2 \pi n q) -
\frac{\pi}{2} \sum_{n=1}^{\infty} \frac{\sin( 2 \pi n q)}{n^{2 m+2} }
\right].
\end{multline}

\medskip

\begin{Lem}
\label{speca}
For $0\le q \le 1$ the function $A_{2}(q)$ is given by
\ba
A_{2}(q) & = & (1- \gamma - \ln 2 \pi) (q^{2} - q + \tfrac{1}{6})
\label{a-two} \\
 &&  - \frac{1}{\pi^{2}}
\sum_{n=1}^{\infty} \frac{ \ln n}{n^{2}} \cos( 2 \pi n q) +
\frac{1}{2\pi}\sum_{n=1}^{\infty} \frac{\sin( 2 \pi n q)}{n^{2} }. \nn
\ea
\no
In particular,
\ba
A_{2}(1) & = & 2\zeta'(-1), \nn \\
A_{2}(\tfrac{1}{2}) & = & - \zeta'(-1) - \tfrac{1}{12} \ln 2, \nn
\\
A_{2}(\tfrac{1}{4}) & = & -\tfrac{1}{4} \zeta'(-1) + \frac{G}{2 \pi},  \nn
\ea
where $G$ is Catalan's constant.
\end{Lem}
\begin{proof}
The expression (\ref{a-two}) follows directly from (\ref{ameven}) at $m=0$.
\end{proof}

Special values of the derivative of the Hurwitz zeta function
at the rational arguments
$q=\tfrac{1}{2},\tfrac{2}{3},\tfrac{1}{4},\tfrac{3}{4},\tfrac{1}{6}$ and
$\tfrac{5}{6}$, for $k$ odd, have been given in \cite{milleradam}.
Special values of $A_k(q)$ for $q>1$ can be obtained by using
Lemma \ref{recu2} a sufficient number of times.

\begin{Example}
\ba
A_{2}(2) & = & 2\zeta'(-1), \nn \\
A_{2}(3) & = & 2\zeta'(-1) + 4\ln 2,\nn \\
A_{2}(\tfrac{3}{2}) & = & - \zeta'(-1) - \tfrac{13}{12} \ln 2, \nn
\\
A_{2}(\tfrac{5}{4}) & = & -\tfrac{1}{4} \zeta'(-1) + \frac{G}{2 \pi}
- \ln 2.  \nn
\ea
\end{Example}

Several of the indefinite integrals derived in the next sections
can be conveniently expressed in terms of a family of functions closely
related to the {\em negapolygamma} family introduced by Gosper
\cite{gosper} as
\begin{equation}
\begin{split}
\psi_{-1}(q) & := \ln \Gamma(q), \\
\psi_{-k}(q) & := \int_0^q \psi_{-k+1}(t) dt,\quad k\ge 2.\\
\stboth
\end{split}
\end{equation}
These functions were later reconsidered by Adamchik \cite{adam1} in the form
\ba\label{adam-negapoly}
\psi_{-k}(q) = \frac{1}{(k-2)!}\int_0^q (q-t)^{k-2}\ln\Gamma(t) dt,\quad k\ge
2.
\ea
\no
We introduce the {\em balanced negapolygamma} functions by
\\
\ba\label{bal-negapolygamma}
\psi ^{( -k)} (q): = \frac{1}{k!}\left[{A_k (q)} - H_{k-1} B_k(q)\right],
\ea
\\
where $k\in\allN$ and $H_r$ is the harmonic number ($H_0:=0$). For
instance,
\\
\ba\label{negapolygamma-minusone}
\psi^{(-1)}(q)&=& A_1(q)=\ln\Gamma(q)+\zeta'(0),\\
\psi^{(-2)}(q)&=& \frac{1}{2}A_2(q)-\frac{1}{2}B_2(q).
\ea
{\bf Note.} We can express the derivative of the Hurwitz zeta
function at non-positive integers as
\ba
\zeta'(-r,q)=r!\psi^{(-1-r)}(q)+\frac{H_r}{1+r} B_{1+r}(q),\quad r\in\nnN.
\ea
\medskip

\begin{Lem}
The balanced negapolygammas can be expressed as
\ba\label{negapolygamma-alt-form}
\psi^{(-k)}(q)  =  e^{- \gamma z } \frac{\partial}{\partial z}
\left[ e^{\gamma z} \frac{\zeta(z,q) }{\Gamma(1-z)} \right]
\Bigg{|}_{z=1-k}.
\ea
\end{Lem}
\begin{proof}
Perform the derivative and use $\psi(k)=H_{k-1}-\gamma$, in
addition to the definition (\ref{defAk}) and the identity
(\ref{zetaber}).
\end{proof}

\no
We shall now study some of the properties of the balanced
negapolygamma functions $\psi ^{( -k)} (q)$. The adjective {\em
balanced} is motivated by the following result:

\begin{Lem}
For $k\in\allN$,
\ba\label{negapolygamma-int-0-1}
\ione\psi^{(-k)}(q)\,dq=0.
\ea
\end{Lem}
\begin{proof}
Use (\ref{Ak-int-0-1}) and the analogous result for the Bernoulli
polynomials,
\ba\nn
\ione B_k(q)\,dq=0,
\ea
for $k\in\allN$.
\end{proof}

\begin{Lem}
For $k\in\allN$,
\ba\label{der-negapolygamma}
\frac{d}{dq}\psi^{(-k)}(q)=\psi^{(-k+1)}(q).
\ea
\end{Lem}
\begin{proof}
For $k\ge 2$, (\ref{der-negapolygamma}) is a direct consequence of
Lemma \ref{der-Ak} and the well-known property of the Bernoulli
polynomials,
\ba
\frac{d}{dq}\,B_k(q)=k B_{k-1}(q).
\ea
The case $k=1$ follows directly from
(\ref{negapolygamma-minusone}) since
\ba\nn
\frac{d}{dq}\ln\Gamma(q)=\psi(q)=\psi^{(0)}(q).
\ea
\end{proof}

The functions $\psi^{(-k)}(q)$ defined by
(\ref{bal-negapolygamma}) are thus
closely related to the ones introduced by Gosper.
In fact, the precise relationship between both families of
polygammas of negative order is
\ba
\psi^{(-k)}(q) =\psi_{-k}(q)+
\sum_{r=0}^{k-1}\frac{q^{k-1-r}}{r!(k-1-r)!}
\left[\zeta'(-r)+H_r\zeta(-r)\right],
\stoli
\ea
where we have used the evaluation
\ba\label{negapolygamma-at-zero}
\psi^{(-1-r)}(0)
=\frac{1}{r!}\left[\zeta'(-r)+H_r\zeta(-r)\right],
\qquad r\in\allN.
\ea
\no
{\bf Note.} According to \cite{adam1},
\ba\label{glaisher-constants}
\zeta'(-r)+H_r\zeta(-r)=-\ln A_r,
\ea
where the $A_r$ are the generalized Glaisher constants defined by
Bendersky \cite{bendersky}.
\\

\begin{Lem}\label{Fourier-negapolygamma}
The balanced negapolygamma functions $\psi^{(-k)}(q)$ admit the
following Fourier expansion in the range $0\le q \le 1$:
\begin{multline}\label{eq-Fourier-negapolygamma}
\psi^{(-k)}(q)=
\frac{2}{(2\pi )^k }\left[ \sum\limits_{n = 1}^\infty
\frac{\ln (2\pi n) + \gamma }{n^k }\cos (2\pi nq -
\frac{k\pi}{2})\right.\\
\left. - \frac{\pi }{2}\sum\limits_{n = 1}^\infty
{\frac{1}{{n^k }}\sin (2\pi nq - \frac{{k\pi }}{2})}  \right].
\end{multline}
\end{Lem}
\begin{proof}
In the definition (\ref{bal-negapolygamma}) of $\psi^{(-k)}(q)$
substitute, depending on the parity of $k$, the Fourier expansions
(\ref{amodd}) or (\ref{ameven}) for the functions $A_k(q)$ and
(\ref{fourier-berodd}) or (\ref{fourier-bereven}) for the
Bernoulli polynomials.
\end{proof}

\begin{Lem}
For $k\in\allN$,
\ba\label{shift-negapolygamma}
\psi^{(-k)}(q+1)=\psi^{(-k)}(q)+\frac{q^{k-1}}{(k-1)!}[\ln
q-H_{k-1}].
\ea
\end{Lem}
\begin{proof}
In the definition (\ref{bal-negapolygamma}) use (\ref{recu2}) and
the property
\ba\nn
B_{m}(q+1)=B_{m}(q)+m\,q^{m-1},
\ea
satisfied by the Bernoulli polynomials.
\end{proof}
\begin{Cor}\label{negapolygamma-at-one}
For $r\in\allN$,
\ba
\psi^{(-1-r)}(1)=\psi^{(-1-r)}(0)
=\frac{1}{r!}\left[\zeta'(-r)+H_r\zeta(-r)\right].
\ea
\end{Cor}
\begin{proof}
Set $k=1+r$ and evaluate (\ref{shift-negapolygamma}) at $q=0$.
Then use (\ref{negapolygamma-at-zero}).
\end{proof}

\section{Indefinite integrals of polygamma functions}
\label{polygamma}

Integral formulae involving the polygamma functions can be
obtained from the corresponding ones for $\zeta(z,q)$, like
(\ref{mom1}) or (\ref{momber1}), by taking the limit $z\to
m=2,3,\ldots$. This limit is not trivial in general, because of
the vanishing of the Pochhammer symbol $(1-z)_{j+1}$ for $j>m-2$,
and the appearance of the function $\zeta(1,q)$ on the right-hand
side of the formulae cited above. Similarly, differentiation at
$z=1-k$ leads to evaluations involving the functions $A_k(q)$ and
the negapolygamma functions. Because of the connection
(\ref{negapolygamma-minusone}) between $\psi^{(-1)}(q)$ and
$\ln\Gamma(q)$, and of the latter with the function $\ln\sin(\pi
q)$, we also obtain indefinite integral formulae involving these
last two functions.
\medskip

The next theorem gives the moments of the polygamma functions
in terms of themselves and the balanced negapolygammas defined
in Section \ref{speciala}. The proof of the theorem rests on the
following result:

\begin{Lem}\label{ber-identity}
Let $p\in\nnN$. Then
\ba
\sum\limits_{r = 0}^p {( - 1)^r \binom{{p + 1}}{{r + 1}}\frac{{q^{p
- r} B_{r + 1} (a + bq)}}{{b^{r + 1} }} = q^{p + 1}  + ( - 1)^p
\frac{{B_{p + 1} (a)}}{{b^{p + 1} }}}.
\ea
\end{Lem}
\begin{proof}
The basic identity
\ba
B_{p+1}(x+y) & = & \sum_{r=0}^{p+1} \binom{p+1}{r} B_{r}(x)y^{p+1-r}
\label{classber}
\ea
\no
appears in \cite{atlas} 19:5:3. The identity of the lemma is obtained by
replacing $x$ by $a+bq$ and $y$ by $-bq$ in (\ref{classber}).
\end{proof}

\begin{Thm}\label{polygamma-moments-thm}
Let $n\in\nnN$, $m \in \mathbb{N}$. Then
\ba
\label{polygamma-moments}
\int q^{n} \psi^{(m)}(a+bq)\,dq = n!
\sum_{j=0}^{n} \frac{(-1)^j}{b^{j+1}(n-j)!} q^{n-j}\psi^{(m-j-1)}(a+bq).
\ea
\end{Thm}
\begin{proof}
In view of relationship (\ref{polygamma-def}),
we set $z=m+\epsilon$ in (\ref{mom1}) and take the limit
$\epsilon\to 0$. For a general value of $n$ there will be three
kinds of terms to consider in the sum on the right-hand side of (\ref{mom1}):

\no
a) For $0\le j\le m-2$ no singularities will arise and we
simply have
\[
\lim _{\epsilon  \to 0} \frac{{\zeta (m - j+\epsilon,q)}}{{( - m -
\epsilon )_{j + 1} }} = \frac{{( - 1)^{m + 1} }}{{m!}}\psi ^{(m
- j - 1)} (q).
\]
b) For $j=m-1$ we use the results
\[
( - m - \epsilon )_m  = ( - 1)^m m!\left[ {1 + H_m \epsilon  + O(\epsilon ^2 )} \right]
\]
\no
and
\[
\zeta (1 + \epsilon ,q) = \frac{1}{\epsilon } - \psi (q) + O(\epsilon ^2 ),
\]
where $\psi(q)=\psi^{(0)}(q)$ is the usual digamma function, to
obtain
\[
\left. {\frac{{\zeta (m - j + \epsilon ,q)}}{{( - m -
\epsilon )_{j + 1} }}} \right|_{j = m - 1}  = \frac{{( - 1)^{m
+ 1} }}{{m!}}\left[ { - \frac{1}{\epsilon } + \psi (q) + H_m  +
O(\epsilon )} \right].
\]
c) For $j\ge m$, say $j=m+r$ with $r\ge 0$, we use
\[
( - m - \epsilon )_{m + 1 + r}  = ( - 1)^{m + 1}
m!r!\epsilon \left[ {1 + (H_m  - H_r )\epsilon  +
O(\epsilon ^2 )} \right]
\]
\no
and
\[
\zeta ( - r + \epsilon ,q) = \zeta ( - r,q) +
\epsilon \zeta'( - r,q) + O(\epsilon ^2 )
\]
to obtain
\begin{multline}\nn
\left. \frac{\zeta (m - j + \epsilon ,q)}{( - m -
\epsilon )_{j + 1}} \right|_{j = m + r}  = \frac{{( - 1)^{m
+ 1} }}{{m!}}\bigg[  - \frac{1}{\epsilon }\frac{{B_{r + 1}
(q)}}{{(r + 1)!}} +H_m\frac{{B_{r + 1} (q)}}{{(r + 1)!}}\\
+ \psi^{(-1-r)}(q) + O(\epsilon ) \bigg],
\end{multline}
where we have used the definition (\ref{bal-negapolygamma}) of the
balanced egapolygamma function and (\ref{zetaber}).
When all the terms are added up, we find, in view of Lemma \ref{ber-identity},
that the coefficients of the $1/\epsilon$
singularity and the term proportional to the harmonic number
$H_m$ reduce to a $q$-independent constant, which can be dropped.
This proves the theorem.
\end{proof}

A similar result holds for the moments of the digamma function:
\begin{Thm}
Let $n \in \nnN$. Then
\ba
\label{digamma-moments}
\int q^{n} \psi(a+bq)\,dq = n!
\sum_{j=0}^{n} \frac{(-1)^j}{b^{j+1}(n-j)!} q^{n-j}\psi^{(-j-1)}(a+bq).
\ea
\end{Thm}
\begin{proof}
Use (\ref{digamma-hurwitz}) in (\ref{mom1}) and proceed along the
same lines as the proof above.
\end{proof}
{\bf Note:} Result (\ref{digamma-moments}) can be identified as an
extension of Theorem \ref{polygamma-moments-thm} to the case
$m=0$.
\medskip

\no
{\bf Note:} Adamchik \cite{adam1,adam2} has provided the alternative
representations
\begin{multline}
\int_{0}^{z} x^{n} \psi(x) \, dx = (-1)^{n}
\left( \frac{B_{n+1} H_{n}}{n+1} - \zeta'(-n) \right) \\
+\sum_{k=0}^{n} (-1)^{k} \binom{n}{k}z^{n-k}
\left( \zeta'(-k,z) - \frac{B_{k+1}(z) H_{k}}{k+1} \right),  \nn
\end{multline}
\no
and
\begin{multline}
\int_{0}^{z} x^{n} \psi(x) \, dx =
\sum_{k=0}^{n-1} (-1)^{k} z^{n-k} \binom{n}{k}
\left( \zeta'(-k) - \frac{B_{k+1}(z) H_{k} }{k+1} \right)\\
- \sum_{k=1}^{n} (-1)^{k} k!
\stnk{n}{k} \; \log G_{k+1}(z+1) + (-1)^{n} H_{n}
\frac{B_{n+1}-B_{n+1}(z)}{n+1},
\nn
\end{multline}
\no
where $\textstyle\stnk{n}{k}$ are the Stirling numbers of the
second kind and $G_{k}(z)$ is the multiple Barnes function. The
first representation can be directly obtained from
(\ref{digamma-moments}) at $a=0$ and $b=1$, by explicit evaluation
of the integral between $q=0$ and $q=z$, making use of
(\ref{negapolygamma-at-zero}) for $\psi^{(-j-1)}(0)$ and of the
definition (\ref{bal-negapolygamma}) of the negapolygamma
functions.

\medskip

By differentiating formula (\ref{mom1}) at $z=1-m\,$ we can arrive
at the following result for the moments of the balanced
negapolygamma functions, which extends Theorem
\ref{polygamma-moments-thm} to negative values of $m$:

\begin{Thm}\label{negapolygamma-moments-thm}
Let $n\in\nnN$, $m \in \allN$. Then
\ba
\label{negapolygamma-moments}
\int q^{n} \psi^{(-m)}(a+bq)\,dq = n!
\sum_{j=0}^{n} \frac{(-1)^j}{b^{j+1}(n-j)!} q^{n-j}\psi^{(-m-j-1)}(a+bq).
\ea
\end{Thm}
\begin{proof}
Differentiate (\ref{mom1}) at $z=1-m$ and use the result
\[
\left. {\frac{d}{{dz}}(1 - z)_{j + 1} } \right|_{z = 1 - m}  =  -
(m)_{j + 1} \left[ {H_{m + j}  - H_{m - 1} } \right]
\]
to obtain the following result for the moments of the function
$A_m(q)$:
\begin{multline}
\label{ak-moments}
\int q^{n} A_m(a+bq)\,dq = m!\,n!
\sum_{j=0}^{n} \frac{(-1)^j}{b^{j+1}(n-j)!(m+j+1)!} q^{n-j}
\Big[ {A_{m + j + 1} (a + bq)}\\
{ - \left( {H_{m + j}  - H_{m - 1} } \right)B_{m + j + 1} (a + bq)}
\Big].
\end{multline}
This result is equivalent to the statement of the theorem in view of
definition (\ref{bal-negapolygamma}) of the balanced
negapolygammas and result (\ref{mom3}) for the moments of the
Bernoulli polynomials.
\end{proof}

As a particular case of the last theorem we obtain a formula for
the moments of the loggamma function $ \ln \Gamma(q)$.
\\

\begin{Example}
\label{exloggamma}
Let $n \in \nnN$ and $a, \, b \in \mathbb{R}$. Then
\begin{multline}
\label{loggamma1}
\int q^{n} \ln \Gamma(a+bq) \, dq  =
\lp \frac{q^{n+1}}{n+1}
+n! \sum_{j=0}^{n} \frac{(-1)^{j} q^{n-j}}{b^{j+1}(n-j)!}
\psi^{(-2-j)}(a+bq).\\
\end{multline}
This generalizes Gosper's result \cite{gosper}, which
establishes that all integrals of the form $\int q^n \ln q! \, dq$
are expressible in terms of $\zeta'(-j,q)$, with $1\le j \le n+1$.
\medskip

\no
As particular cases of (\ref{loggamma1}) we have
\begin{multline}
\int q^{n} \ln \Gamma(q) \, dq = \lp \frac{q^{n+1}}{n+1}
+n! \sum_{j=0}^{n} \frac{(-1)^{j} q^{n-j}}{(n-j)!}
\psi^{(-2-j)}(q) \\
\label{loggamma2}
\end{multline}
\no
and
\begin{multline}
\int q^{n} \ln \Gamma(1-q) \, dq = \lp \frac{q^{n+1}}{n+1}
- n!\sum_{j=0}^{n} \frac{q^{n-j}}{(n-j)!}
\psi^{(-2-j)}(1-q).\\
\label{loggamma3}
\end{multline}
\end{Example}
\begin{proof}
Apply theorem (\ref{negapolygamma-moments-thm}) to
$\psi^{(-1)}(q)=\ln\Gamma(q)-\lp$.
The expressions
(\ref{loggamma2}) and (\ref{loggamma3})
correspond to $a=0, \, b=1$ and $a=1, \, b=-1$ respectively.
\end{proof}

The two special cases of Example \ref{exloggamma} are now combined with
the reflection formula for the gamma function
\ba
\Gamma(q) \Gamma(1-q) & = & \frac{\pi}{\sin \pi q} \label{reflec}
\ea
\no
to obtain an expression for the moments of $ \ln \sin \pi q$. \\

\begin{Example}
Let $n \in \nnN$. Then
\begin{multline}
\label{logsine1}
\int q^{n} \ln \sin \pi q \; dq = -\frac{q^{n+1} \ln 2}{n+1}\\
- n!\sum_{j=0}^{n} \frac{q^{n-j}}{(j+2)!(n-j)!} \left[ (-1)^{j} A_{j+2}(q) -
A_{j+2}(1-q) \right].
\end{multline}
\end{Example}
\begin{proof}
Use the reflection formula for $\Gamma(q)$, results
(\ref{loggamma2}, \ref{loggamma3}) and the definition
(\ref{bal-negapolygamma}) to produce (\ref{logsine1}). The term
that corresponds to the Bernoulli polynomials disappears in view
of
\ba
(-1)^{j} B_{j+2}(q) & = & B_{j+2}(1-q). \label{vanish}
\ea
\end{proof}

\medskip

\begin{Example}
\ba\label{logsine2}
\int e^{q} \ln \sin \pi q \; dq  =  - e^{q} \left[ \ln 2 +
 \sum_{j=0}^{\infty} \frac{(-1)^{j} A_{j+2}(q) - A_{j+2}(1-q)}{(j+2)!} \right].
\ea
\no
\end{Example}
\begin{proof}
Divide (\ref{logsine1}) by $n!$ and sum over $n$.
\end{proof}

\medskip

\no
\begin{Example}
Integrating (\ref{logsine2}) by parts yields
\begin{multline}\label{cotg-1}
\int e^{q} \; \text{cotg} \, \pi q \; dq = \frac{e^{q}}{\pi}
\left[ \ln \sin \pi q + \ln 2 +
 \sum_{j=0}^{\infty} \frac{(-1)^{j} A_{j+2}(q) - A_{j+2}(1-q)}{(j+2)!} \right].
\end{multline}
\end{Example}

\section{Some definite integrals}
\label{definite}

Some of the definite integral formulae given in \cite{esmo1}, in the
range $(0,1)$, can be obtained directly from the indefinite
integral formulae given in Sections \ref{indefinite} and \ref{polygamma}.

\begin{Example}
Evaluate equation (\ref{mom1}) between 0 and 1 to obtain
\begin{multline}
\label{mom1-def}
\ione q^{n} \zeta(z,a+bq) \, dq =
n! \sum_{j=0}^{n-1} \frac{(-1)^j \zeta(z-j-1,a+b)}{b^{j+1}(1-z)_{j+1}
(n-j)!}\\
+\frac{n!(-1)^n}{b^{n+1}(1-z)_{n+1}}\left(\zeta(z-n-1,a+b)-\zeta(z-n-1,a)\right).
\end{multline}
\end{Example}
\no
As a particular case we obtain formula (12.2) of
\cite{esmo1}:
\begin{Cor}
Let $n\in\nnN$ and $z\in\allR$, with $z-n-1<0$. Then
\ba
\ione q^{n} \zeta(z,q) \, dq =
n! \sum_{j=0}^{n-1} \frac{(-1)^j \zeta(z-j-1)}{(1-z)_{j+1}
(n-j)!}.
\ea
\end{Cor}
\begin{proof}
Set $b=1$ in (\ref{mom1-def}) and use the identity
(\ref{zeta-shift}) and the hypothesis $z-n-1<0$ to get rid of the
last term in (\ref{mom1-def}) in the limit $a\to 0$.
\end{proof}

The evaluation of formula (\ref{logsine1}) between $q=0$ and $q=1$
leads to formula (5.6) of \cite{esmo1}:
\begin{Example}
Let $n \in \nnN$. Then
\ba
\label{momlogsine-def}
\ione q^{n} \ln (\sin \pi q) dq = - \frac{\ln 2}{n+1}
+ n! \sum_{k=1}^{\lfloor{ \tfrac{n}{2} \rfloor} }
\frac{(-1)^{k} \zeta(2k+1)}{(2 \pi)^{2k} \, (n+1 - 2k)!}.
\ea
\end{Example}
\begin{proof}
Direct evaluation of the right-hand side of (\ref{logsine1})
gives, in view of property (\ref{Ak-at-zero}),
\begin{multline}
\ione q^{n} \ln \sin \pi q \; dq = -\frac{\ln 2}{n+1}\\
- n!\sum_{j=0}^{n-1} \frac{1}{(j+1)!(n-j)!} \left[ (-1)^{j}-1 \right]
\zeta'(-j-1).
\end{multline}
Clearly only the terms with $j$ odd, say $j=2k-1$, survive in the sum.
(\ref{momlogsine-def}) then follows directly from (\ref{berzeta2}).
\end{proof}

\begin{Example}
Let $n \in \nnN$. Then
\begin{multline}
\label{halfmomlogsine-def}
\ionehalf q^{n} \ln (\sin \pi q) dq = -\frac{1}{2^{n+1}}
\Biggl[\frac{\ln 2}{n+1}\\
+ n! \sum_{k=1}^{\lfloor{ \tfrac{n+1}{2} \rfloor} }
\frac{(-1)^{k}(2^{2k}-1) \zeta(2k+1)}{(2 \pi)^{2k} \, (n+1 - 2k)!}\Biggr]
-\frac{1-(-1)^n}{n+1}\zeta'(-n-1).
\end{multline}
\end{Example}
\begin{proof}
Use (\ref{Ak-at-onehalf}) in (\ref{logsine1}).
\end{proof}
\no
For instance,
\ba
\ionehalf q \ln (\sin \pi q) dq &=&
 - \frac{1} {8}\ln 2 + \frac{{7\zeta (3)}} {{16\pi ^2 }},\nn\\
\ionehalf q^{2} \ln (\sin \pi q) dq &=&
 - \frac{1} {{24}}\ln 2 + \frac{{3\zeta (3)}} {{16\pi ^2 }},\nn\\
\ionehalf q^{3} \ln (\sin \pi q) dq &=&
 - \frac{1} {{64}}\ln 2 + \frac{{9\zeta (3)}}
{{64\pi ^2 }} - \frac{{93\zeta (5)}} {{128\pi ^4 }}.\nn
\ea

\no
\begin{Example}
Formulae (\ref{loggamma1}) or (\ref{loggamma2}) allow us to derive
Gosper's formulae for integrals of $\ln\Gamma(q)$ \cite{gosper} in
a very economical way.
For instance, setting $n=0$, $a=b=1$ in (\ref{loggamma1}) yields
\begin{multline}
\label{loggamma1-def}
\int_{0}^{q} \ln \Gamma(q+1) \, dq  =
q\lp + \frac{1}{2} A_{2}(q+1) - \frac{1}{2} B_{2}(q+1)
-\zeta'(-1) + \frac{1}{2} B_{2}.
\end{multline}
Evaluation of the right-hand side at $q=\tfrac{1}{2}$ and
$\tfrac{1}{4}$ yields, respectively,
\ba
\int_{0}^{1/2} \ln \Gamma(q+1) dq &=&- \frac{3}{8} - \frac{{13}}{{24}}\ln
2 + \frac{{1}}{2}\lp - \frac{3}{2}\zeta '( - 1),
\stboth\nn
\\
\int_{0}^{1/4} \ln \Gamma(q+1) dq &=&- \frac{5}{32} - \frac{{1}}{{2}}\ln
2 + \frac{{1}}{4}\lp - \frac{9}{8}\zeta '( - 1)+\frac{G}{4\pi},\nn
\ea
which can easily be seen to be equivalent to Gosper's formulae, after
using Riemann's functional equation for the Riemann zeta function to
express $\zeta'(-1)$ in the form
\ba\label{zetaprimeofone}
\zeta '( - 1) = \frac{{\zeta '(2)}}{{2\pi ^2 }} -
\frac{1}{{12}}(2\lp +\gamma  - 1).
\ea
\end{Example}

\medskip

\no
\begin{Example}
For $k,k'\in\allN$,
\begin{multline}\nn
\ione\psi ^{( - k)} (q)\psi ^{( - k')} (q)\,dq =
\frac{{2\cos (k - k')\frac{\pi}{2}}}{{(2\pi )^{k + k'} }}
\Bigg[ \zeta'' (k + k') - 2(\gamma  +
\ln 2\pi )\zeta' (k + k')\\
 + \left\{ {(\gamma  + \ln 2\pi )^2  +
\frac{{\pi ^2 }}{4}} \right\}\zeta (k + k') \Bigg].
\end{multline}
\no
The special case $k = k' =1$ reduces to
\ba
\int_{0}^{1} \left( \ln \Gamma(q) \right)^{2} dq & = &
\frac{\gamma^{2}}{12} + \frac{\pi^{2}}{48} + \frac{1}{3} \gamma
\ln \sqrt{2 \pi} + \frac{4}{3} \ln^{2} \sqrt{2 \pi} \nn \\
&&\phantom{xxx}-
(\gamma + 2 \ln \sqrt{2 \pi} ) \frac{\zeta'(2)}{\pi^{2}} +
\frac{\zeta''(2) }{2 \pi^{2}}, \nn
\ea
\no
given in \cite{esmo1}.
\begin{proof}
Use the representation (\ref{negapolygamma-alt-form}) of the balanced
negapolygammas to obtain the
desired result, by direct differentiation of the formula
\begin{multline}
\ione \zeta(z',q) \zeta(z,q) dq  = \frac{2 \Gamma(1-z) \Gamma(1-z')}
{(2 \pi)^{2 - z - z'}} \zeta(2 - z - z')
\cos \left( \frac{\pi(z-z')}{2} \right),
\end{multline}
valid for real $z, \, z'\le 0$, given in \cite{esmo1}.
\end{proof}
\end{Example}

\no
{\bf Acknowledgments}. The authors would like to thanks G. Boros for many
suggestions. The first author would like to thank the Department of
Mathematics at Tulane University for its hospitality and the support of CONICYT
(Chile) under grant P.L.C. 8000017. The second author acknowledges the
partial support of NSF-DMS 0070567, Project number 540623.

\end{document}